\documentclass[10pt]{amsart}
\usepackage{amsmath, amsthm, amssymb}
\usepackage{graphicx}
\usepackage{mathbbol}
\usepackage{txfonts}
\usepackage[usenames]{color}

\newtheorem{thm}{Theorem}[section]
\newcommand{\bt}{\begin{thm}}
\newcommand{\et}{\end{thm}}

\newtheorem{cor}[thm]{Corollary}   
\newcommand{\bc}{\begin{cor}}
\newcommand{\ec}{\end{cor}}

\newtheorem{lem}[thm]{Lemma}   
\newcommand{\bl}{\begin{lem}}
\newcommand{\el}{\end{lem}}

\newtheorem{prop}[thm]{Proposition}
\newcommand{\bp}{\begin{prop}}
\newcommand{\ep}{\end{prop}}

\newtheorem{defn}[thm]{Definition}
\newcommand{\bd}{\begin{defn}}    
\newcommand{\ed}{\end{defn}}

\newtheorem{rmrk}[thm]{Remark}   

\newcommand{\br}{\begin{rmrk}}
\newcommand{\er}{\end{rmrk}}

\newtheorem{example}[thm]{Example}

\newcommand{\be}{\begin{equation}}

 \newcommand{\ee}{\end{equation}}

\newcommand{\N}{\mathbb{N}}

\newcommand{\E}{\mathbb{E}}

\newcommand{\vare}{\varepsilon}

\newcommand{\diam}{\operatorname{Diam}}

\newcommand{\Slice}{\operatorname{Slice}}

\newcommand{\disjointunion}{\sqcup}

\newcommand{\SF}{{\mathbf {SF}}}

\newcommand{\vol}{\operatorname{Vol}}

\newcommand{\fillvol}{{\operatorname{FillVol}}}

\newcommand{\spt}{\operatorname{spt}}

\begin{document}

\title{The Tetrahedral Property and a new Gromov-Hausdorff
Compactness Theorem}

\author{Christina Sormani}
\thanks{partially supported by NSF DMS 1006059 and a PSC CUNY Research Grant.}
\address{CUNY Graduate Center and Lehman College}
\email{sormanic@member.ams.org}




\begin{abstract}
We present the Tetrahedral Compactness Theorem
which states that sequences of Riemannian manifolds with a uniform
upper bound on volume and diameter that satisfy a uniform
tetrahedral property have a subsequence which converges in
the Gromov-Hausdorff sense to a countably $\mathcal{H}^m$
rectifiable metric space of the same dimension.   The
tetrahedral property depends only on distances between points
in spheres, yet we show it provides a lower bound on the
volumes of balls.  The proof is based upon intrinsic flat convergence
and a new notion called the sliced filling volume of a ball.
\end{abstract}

\maketitle

\section{Introduction}

\begin{figure}[htbp] 
   \centering
   \includegraphics[width=3in]{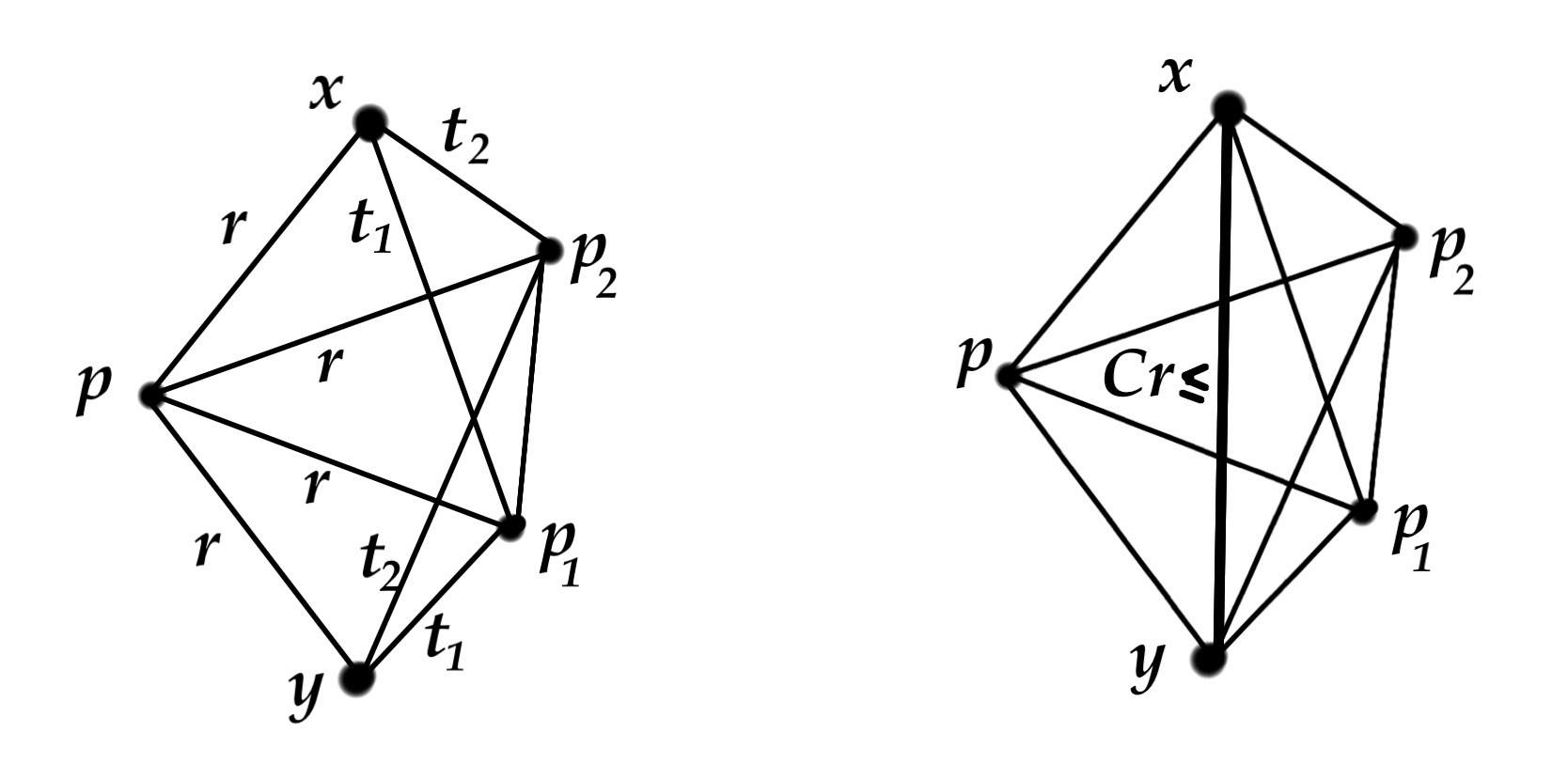} 
   \caption{Tetrahedral Property in 3D}
   \label{fig-tetra-prop}
\end{figure}

We introduce the tetrahedral property 
which is an estimate on tetrahedra (see Figure~\ref{fig-tetra-prop}):

\vspace{.3cm}
\begin{defn} \label{defn-tetra}
Given $C>0$ and $\beta\in (0,1)$, a metric space $X$
has the $m$ dimensional
$C,\beta$-tetrahedral property at a point $p$
for radius $r$ if
one can find points $p_1,...p_{m-1}\subset \partial B_p(r)\subset\bar{X}$,
such that 
\be
h(p,r, t_1,...,t_{m-1}) \ge Cr \qquad \forall (t_1,...,t_{m-1}) \in [(1-\beta)r, (1+\beta)r]^m
\ee
where 
$
h(p,r, t_1,..., t_{m-1})=
\inf\left\{d(x,y): \,\, x\neq y,\,\, x,y\in  P(p, r, t_1,...,t_{m-1})\right\}
$
when
\be
P(p,r, t_1,..., t_{m-1})= \rho_p^{-1}(r) \cap \rho_{p_1}^{-1}(t_1)\cap
\cdots \cap \rho_{p_{m-1}}^{-1}(t_{m-1}) \neq \emptyset
\ee
and $h(p,r, t_1,..., t_{m-1})=0$ otherwise.
In particular $P(p,r,t_1,..., t_{m-1})$ is a discrete set of
points.
\end{defn}

As this property is rather strong, we introduce the
integral tetrahedral property:

\vspace{.3cm}
\begin{defn} \label{defn-int-tetra}
Given $C>0$ and $\beta\in (0,1)$, a metric space $X$
is said to have the $m$ dimensional integral
$C,\beta$-tetrahedral property at a point $p$
for radius $r$ if
$\exists \,\,p_1,...p_{m-1}\subset \partial B_p(r)\subset\bar{X}$,
such that 
\be
\int_{t_1=(1-\beta)r}^{(1+\beta)r} \cdots \int_{t_{m-1}=(1-\beta)r}^{(1+\beta)r }
h(p,r,t_1,...t_{m-1}))
\, dt_1dt_2...dt_{m-1}
\ge C (2\beta)^{m-1}r^m.
\ee
\end{defn}

\vspace{.3cm}

We can prove that both of these properties
provide an estimate on volume:

\vspace{.3cm}
\begin{thm}\label{tetra-manifold}
If $p_0$ lies in a Riemannian manifold 
that has the
$m$ dimensional (integral)
$C,\beta$-tetrahedral property at a point $p$
for radius $R$ 
then
\be
\vol(B(p,r))\ge C (2\beta)^{m-1}r^m
\ee
\end{thm}

\vspace{.3cm}
As a consequence of Gromov's Compactness Theorem we then
have:

\vspace{.3cm}
\begin{thm}\label{tetra-compactness}
Given $r_0>0, \beta\in (0,1), C>0, V_0>0$.  If
a sequence of compact Riemannian manifolds, $M^m$, 
has $\vol(M^m) \le V_0$,
$\diam(M^m) \le D_0$, 
and the $C, \beta$ (integral) tetrahedral
property for all balls of radius $\le r_0$, then a subsequence
converges in the Gromov-Hausdorff sense.   In particular
they have a uniform upper bound on diameter
depending only on these constants.   
\end{thm}

\vspace{.3cm}

\begin{rmrk}\label{tetra-compactness-r}
In fact we prove there is an intrinsic flat limit as well
and the intrinsic flat and GH limit agree.  Thus the
limit space in Theorem~\ref{tetra-compactness}
is a countably $\mathcal{H}^m$ rectifiable
metric space.
\end{rmrk}

\section{Examples}

\begin{example} \label{ex-torus-Euclidean}
On Euclidean space, $\E^3$, taking $p_1, p_2 \in \partial B(p,r)$ to
such that $d(p_1,p_2)=r$, then there exists exactly two
points $x,y \in P(p,r,r, r)$ each forming a tetrahedron with
$p, p_1, p_2$.  As we vary $t_1, t_2\in (r/2, 3r/2)$, we still have
exactly two points in $P(p,r,t_1, t_2)$.   By scaling we see that
\be
h(p,r, t_1, t_2)= r h(p,1, t_1/r, t_2/r) \ge 
C_{\E^3} r
\ee
where
$
C_{\E^3}=\inf\{h(p,1,s_1, s_2): \, s_i\in (1/2, 3/2)\}>0.
$
Taking $\beta=1/2$,  we see
that $\E^3$ satisfies the $C_{\E^3}, \beta$ tetrahedral property.
\end{example}

\vspace{.3cm}
\begin{example} \label{ex-torus-tetra}
On a torus, $M_\epsilon^3=S^1 \times S^1 \times S_\epsilon^1$ where
$S_\epsilon^1$ has been scaled to have diameter $\epsilon$
instead of $\pi$, we see that $M^3$ satisfies the
$C_{\E^3}, (1/2)$ tetrahedral property at $p$ for all $r< \epsilon/4$.
By taking $r<\epsilon/4$, we guarantee that the shortest
paths between $x$ and $y$ stay within the ball $B(p,r)$ allowing
us to use the Euclidean estimates.   If $r$ is too large,
$P(p, r, t_1, t_2)=\emptyset$.
So for a sequence $M_\epsilon$ with $\epsilon \to0$ we fail to
have a uniform tetrahedral property.  There is a Gromov-Hausdorff
limit but it is not three dimensional.
\end{example}

\vspace{.3cm}
\begin{example}
Suppose one creates a Riemannian manifold $M_\vare^3$, by gluing
together two copies of Euclidean space with a large collection
of tiny necks between corresponding points.  That is,
\be
M_\vare^3= \left(\mathbb{E}^3 \setminus \bigcup B_{z_i}(\vare)\right)   
\disjointunion \left(\mathbb{E}^3 \setminus \bigcup B_{z_i}(\vare)\right)   
\ee
where points on $\partial B_{z_i}(\vare)$ in the first copy
of Euclidean space are joined to corresponding points on
$\partial B_{z_i}(\vare_i)$ in the second copy
of Euclidean space.  We choose $z_i$ such that
$\mathbb{E}^3 \subset \bigcup_{i=1}^\infty B_{z_i}(10\vare)$
and the balls $B_{z_i}(\vare)$ are pairwise disjoint.
Then for $r>>\vare$, 
we will have an $x$ and a $y$ as in $\mathbb{E}^3$, 
but we will also
have a nearby $x'$ and $y'$ in the second copy, with
$
d(x,x')<20\vare.
$
So for a sequence $M_\vare$ with $\vare\to0$ we fail to
have a uniform tetrahedral property.   If we create $M_\vare^3$
by joining increasingly many copies of Euclidean space together,
this sequence wouldn't even have a subsequence converging
in the Gromov-Hausdorff sense.
\end{example}

\section{Intrinsic Flat Convergence}

The Intrinsic Flat distance between Riemannian manifolds 
was introduced by the author and Stefan Wenger in \cite{SorWen2}.
It was defined using Gromov's idea of
isometrically embedding two Riemannian manifolds into a
common metric space.  Rather than measuring the Hausdorff 
distance between the images as Gromov did
when defining the Gromov-Hausdorff distance in
\cite{Gromov-metric}, one views the images as
integral currents in the sense of Ambrosio-Kirchheim in
\cite{AK} and takes the flat distance between them.
The author and Wenger proved that intrinsic flat
limit spaces are countably $\mathcal{H}^m$ rectifiable
metric spaces in \cite{SorWen2}.

Wenger has proven a compactness theorem for intrinsic
flat convergence \cite{Wenger-Compactness}, but we do 
not need to apply that compactness theorem to prove
the compactness theorems stated here.  Instead our compactness
theorem is based upon the Gromov-Hausdorff compactness
theorem \cite{Gromov-metric}
and the fact that we obtain a uniform lower bound on
the volumes of balls [Theorem~\ref{tetra-manifold}].   
Applying Ambrosio-Kirchheim's Compactness Theorem
of \cite{AK}, Wenger and the author proved that once a sequence of manifolds
converges in the Gromov-Hausdorff sense to a limit space $Y$,
then a subsequence converges in the Intrinsic Flat sense
to a subset, $X$, of $Y$ \cite{SorWen2}.  In \cite{SorWen1},
estimates on the filling volumes of spheres were applied to
prove the two limit spaces were the same when the sequence
of manifolds has nonnegative Ricci curvature.   Recall that
filling volumes were introduced by Gromov in \cite{Gromov-filling}.

Here we do not have strong estimates on the filling volumes 
of spheres.   To prove Theorem~\ref{tetra-manifold} we
first define the sliced filling volumes of balls and then prove a compactness
theorem:

\vspace{.3cm}
\begin{defn}\label{sliced-filling-vol}
Given points $q_1,...,q_k\in M^m$, where $k<m$,
with distance functions $\rho_i(x)=d(x, q_i)$,
we define the sliced filling volume of a sphere
$\partial B(p,r)$, to be
\be
\SF(p,r,q_1,...,q_k)
=\int_{t_1=m_1}^{M_1}\int_{t_2=m_2}^{M_2}\cdots \int_{t_k=m_k}^{M_k} \fillvol(\partial\Slice(B(p,r),\rho_1,...\rho_k,t_1,...,t_k)) \, \mathcal{L}^k
\ee
where $m_i=\min\{\rho_i(x): x\in \bar{B}_p(r)\}$
and $M_i=\max\{\rho_i(x): x\in \bar{B}_p(r)\}$
and where the slice is defined as in Geometric Measure
Theory so that it
is supported on
$
B(p,r) \cap  \rho_1^{-1}(t_1) \cap \cdots \cap \rho_k^{-1}(t_k).
$
\end{defn}

\vspace{.3cm}
\begin{defn} \label{defn-SF_k}
Given $p\in M^m$,
then for almost every $r$, we can define the
$k^{th}$ sliced filling,
\be
\SF_k(p,r) = \sup\big\{ \SF(p, r, q_1,...,q_k): q_i \in \partial B_p(r)\big\}.
\ee
\end{defn}

\vspace{.3cm}
\begin{thm} \label{SF_k-compactness}
Let $V_0, D_0, r_0>0$ and $C(r)>0$.
If $M_i^m$ have  $\vol(M_i) \le V_0$, $\diam(M_i)\le D_0$,
and 
\be
\SF_k(p,r) \ge C(r)>0 \qquad \forall i\in \N, \,\,\forall p\in M_i
\textrm{ and almost every } r\in (0, r_0)
\ee
then a subsequence of the $M_i$ converges in the Gromov-Hausdorff
sense to a limit space which is also the intrinsic flat
limit of the sequence and is thus a countably $\mathcal{H}^m$
rectifiable metric space.
\end{thm}

\vspace{.3cm}
This theorem is proven by the author in \cite{Sormani-Properties}.
We first show that $\vol(B(p,r)) \ge \SF_k(p,r)$, so that a subsequence
has a Gromov-Hausdorff limit, $M_\infty$, 
by Gromov's Compactness Theorem.
We next observe that when the points $p_j \in M_j$ converge to
a point $p_\infty$ in the Gromov-Hausdorff limit, $M_\infty$, 
their sliced fillings converge.   Applying Ambrosio-Kirchheim's
Slicing Theorem, we can then estimate the mass of the limit
current and prove that the Gromov-Hausdorff and Intrinsic Flat
limits agree.  

\section{Tetrahedral Property}

The final step in the proof of the Theorem~\ref{tetra-compactness}
is to relate the tetrahedral property to the $k=m-1$ sliced filling volume.
We first observe that 
\be
\spt \Big(\partial\Slice\big(B(p,r),\rho_1,...\rho_{m-1},t_1,...,t_{m-1})\big)\Big)
\,=\,
\partial B_p(r) \,\cap\, \bigcap_{i=1}^{m-1} \partial B_{q_i}(t_i)
\ee
which is a discrete collection of points for almost every value
of $(t_1,...t_{m-1})$.   So we prove a theorem that the filling volume 
of a $0$ dimensional integral current can be bounded below
by the distance between the closest pair of points in the
current's support.   We then have:

\vspace{.3cm}
\begin{thm}\label{tetra-ball}
If $M^m$ is a Riemannian manifold
with the $m$ dimensional (integral)
$C,\beta$-tetrahedral property at a point $p$
for radius $r$ 
then
$
\vol(B(p,r))\ge \SF_{m-1}(p,r)\ge C (2\beta)^{m-1}r^m
$
\end{thm}

\vspace{.3cm}
Combining this theorem with Theorem~\ref{SF_k-compactness},
we obtain both Theorem~\ref{tetra-compactness} and Remark~\ref{tetra-compactness-r}.  These theorems and related theorems are proven in \cite{Sormani-Properties} which is available on the arxiv.  
That paper will include many additional
results before it is completed, as it explores many properties which 
are continuous under intrinsic flat convergence even in settings where
there are no Gromov-Hausdorff limits and where the spaces are not
Riemannian manifolds.


\end{document}